\documentclass[12pt]{article}

\usepackage{amsmath,amssymb,amsthm}
\usepackage[margin=1.2in]{geometry}

\usepackage[osf]{newtxtext} 
\usepackage{zlmtt}          
\usepackage[libertine,bigdelims,cmintegrals,vvarbb]{newtxmath}
\usepackage[cal=euler,frak=euler]{mathalfa} 

\usepackage{url}                    
\usepackage{tikz}                   

\title{A Quadratic Harmonic Approximation}

\makeatletter
\def\section{\@startsection{section}{1}{\z@}{-3.5ex plus -1ex minus
			  -.2ex}{2.3ex plus .2ex}{\large\bf}}
\def\subsection{\@startsection{subsection}{2}{\z@}{-3.25ex plus -1ex
			  minus -.2ex}{1.5ex plus .2ex}{\normalsize\bf}}
\makeatother

\newcommand{\bull}{{\scriptstyle\bullet}} 
\newcommand{\dsp}{\displaystyle}    
\newcommand{\eps}{\varepsilon}      
\renewcommand{\geq}{\geqslant}      
\renewcommand{\leq}{\leqslant}      
\newcommand{\sst}{\scriptstyle}     

\newtheorem{thm}{Theorem}

\newcommand{\word}[1]{\quad\mbox{#1}\quad} 

\begin{document}

\maketitle

\section{Introduction} 
\label{sec:intro}

Some eight hundred years ago the French archbishop Nicholas Oresme
developed his beautiful proof that the $n$-th \emph{harmonic number}:
\begin{equation*}
\label{eq:hn} 
H_n := 1 + \frac{1}{2} + \frac{1}{3} +\cdots+ \frac{1}{n}
\end{equation*}
satisfies the following growth inequality:
\begin{equation*}
\label{eq:Oresme} 
H_{2^k} > 1 + \frac{k}{2}\,,
\end{equation*}
and thereby presented the first example in the history of mathematics, and the first seen by
 countless generations of calculus students, of an infinite
series $\sum_{n=1}^\infty \frac{1}{n}$ which diverges although its
$n$-th term decreases to zero.

Unfortunately $H_n$ has no (known) simple closed formula
representation and so its further study demanded that mathematicians
find suitable \emph{approximation} formulas. The great Leonhard Euler
applied his famous Euler--Maclaurin sum formula to obtain the following
asymptotic formula:
\begin{equation}
\label{eq:euler} 
H_n \sim \ln n + \gamma + \frac{1}{2n} - \frac{1}{12n^2}
+ \frac{1}{120n^4} - \frac{1}{252n^6} \pm \cdots
\end{equation}
where $\gamma \approx 0.577\dots$ is Euler's constant.  If one truncates this expansion after $n$ terms, then the
\emph{error} $E_n$ one commits in using the truncated series as an approximation to $H_n$ is less than the first term truncated and has the same
sign.

There is considerable interest in proving simplified versions
of~\eqref{eq:euler} without using the heavy analytical machinery
employed by Euler. For example, Robert M. Young~\cite{Young} used an
elegant geometrical argument to prove the \emph{linear} approximation:
\begin{equation*}
\label{eq:young} 
H_n = \ln n + \gamma  + \frac{1}{2(n+\theta_n)}\quad (0 < \theta_n < 1).
\end{equation*}

In this note we will modify his argument to prove the following
\emph{quadratic} approximation.

\begin{thm} 
\label{th:quadratic}
\begin{equation}
\label{eq:quadratic} 
\boxed{H_n = \ln n + \gamma + \frac{1}{2n} - \frac{1}{12n^2} + \eps_n}
\word{where}
0 < \eps_n < \frac{1}{4n^3}\,.
\end{equation}
\end{thm}

Admittedly, the error in Euler's formula satisfies 
$0 < E_4 < 1/120n^4$ which is much sharper; but all known proofs
require much more difficult analysis than ours, while our method still
gives the dominant quadratic term $-1/12n^2$ and so is not too bad.
The interest in our note is the simplicity of method to obtain a rather
difficult result.


\section{Geometrical proof}
\label{sec:geo-proof}

We let $T_n$ be the trapezoid with base the line segment $(n,0)$ to
$(n + 1,0)$ on the $x$-axis, sides the lines $x = n$ and $x = n + 1$
and slanted top the line segment joining the point 
$\bigl( n, \frac{1}{n} \bigr)$ to the point 
$\bigl( n + 1, \frac{1}{n+1} \bigr)$.  We decompose $T_n$ into three parts:
\begin{itemize}
\item
The \emph{rectangle}~$r_n$,  with vertices $\bigl( n, 0),\bigl( n+1, 0 \bigr),\bigl( n+1, \frac{1}{n+1} \bigr),\bigl( n, \frac{1}{n+1} \bigr)$ and area $\frac{1}{n+1}\,$.
\item
The \emph{curvilinear right-angled triangle} with base the top of the rectangle
$r_n$ and side the segment joining $\bigl( n,\frac{1}{n+1} \bigr)$
to $\bigl( n,\frac{1}{n} \bigr)$ and curved ``hypotenuse'' the portion
of the curve $y = \frac{1}{x}$ joining the point
$\bigl( n,\frac{1}{n} \bigr)$ to the point
$\bigl( n + 1,\frac{1}{n+1} \bigr)$. We call its area~$\delta_n$.
\item
The \emph{``sliver''} bounded below by the arc of $y = \frac{1}{x}$ and above
by the top of the trapezoid. We call its area~$\sigma_n$.
\end{itemize}

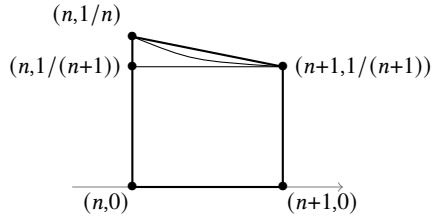
\begin{figure}[htb] 
\centering
\begin{tikzpicture}[x={(0:2cm)}, y={(90:8cm)}]
\coordinate (A) at (4,0) ;
\coordinate (B) at (5,0) ;
\coordinate (C) at (4,0.25) ;
\coordinate (D) at (5,0.2) ;
\coordinate (E) at (4,0.2) ;
\coordinate (F) at (4.4,0.21) ;
\coordinate (X) at (5.4,0) ; \coordinate (Xm) at (3.6,0) ;
\draw[gray, ->] (Xm) -- (X) ; 
\draw[thick] (A) -- (B) -- (D) -- (C) -- cycle ; 
\draw (E) -- (D) ; 
\draw (C) .. controls (F) .. (D) ; 
%
\draw (A) node[below=6pt, left=-3pt]  {$\sst (n,0)$} ;
\draw (B) node[below=6pt, right=-3pt] {$\sst (n+1,0)$} ;
\draw (C) node[above left]  {$\sst (n,1/n)$} ;
\draw (D) node[right]       {$\sst (n+1,1/(n+1))$} ;
\draw (E) node[left]        {$\sst (n,1/(n+1))$} ;
\foreach \pt in {A,B,C,D,E}  \draw (\pt) node {$\bull$} ;
\end{tikzpicture}
\end{figure} 

We define
\begin{equation}
\label{eq:gamma-n} 
\gamma_n := H_n - \ln n.
\end{equation}

Then, as is well known~\cite{Young} (see also~\cite{Havil}),
$$
\sum_{p=n}^\infty \delta_p = \gamma_n - \gamma.
$$
In the interest of completeness we reproduce Young's nice proof:
\begin{align*}
\sum_{p=n}^N \delta_p
&=\left[\int_{n}^{n+1}\frac{1}{x}~dx-\frac{1}{n+1}\right]+\left[\int_{n+1}^{n+2}\frac{1}{x}~dx-\frac{1}{n+2}\right]+\cdots+\left[\int_{N-1}^{N}\frac{1}{x}~dx-\frac{1}{N}\right]\\
&=\int_{n}^{N}\frac{1}{x}~dx-\sum_{r=1}^{N-n}\frac{1}{n+r}=\int_{n}^{N}\frac{1}{x}~dx-\left[\sum_{r=1}^{N}\frac{1}{r}-\sum_{r=1}^{n}\frac{1}{r}\right]\\
&=\left[\ln N-\sum_{r=1}^{N}\frac{1}{r}\right]-\left[\ln n-\sum_{r=1}^{n}\frac{1}{r}\right]
\end{align*}
Now we let $N\rightarrow\infty$ in the last equality and use the definitions of $\gamma_n$ and $\gamma$ to obtain
$$
\sum_{p=n}^\infty \delta_p =-\gamma+\gamma_n=\gamma_n-\gamma
$$
which was to be proved.
But the area of the right-angled triangle at the top of the trapezoid equals
$$
\frac{1}{2} \biggl( \frac{1}{n} - \frac{1}{n+1} \biggr)
= \delta_n + \sigma_n,
$$
and summing from $n$ to infinity we obtain
$$
\frac{1}{2n} = H_n - \ln n - \gamma + \sum_{p=n}^{\infty} \sigma_p,
$$
that is,
$$
H_n = \ln n + \gamma + \frac{1}{2n} -  \sum_{p=n}^{\infty} \sigma_p.
$$
Since $\sigma_n$ is the area of the trapezoid decreased by the area
under the curve $y = 1/x$. we obtain
\begin{align*}
\sigma_n 
&= \frac{1}{2} \biggl( \frac{1}{n} + \frac{1}{n + 1} \biggr)
- \int_n^{n+1} \frac{1}{x} \,dx
= \frac{1}{2n} + \frac{1}{2n(1 + 1/n)} 
- \ln \biggl( 1 + \frac{1}{n} \biggr)
\\
&= \frac{1}{2n} + \left[\frac{1}{2n} - \frac{1}{2n^2} + \frac{1}{2n^3}
- \frac{1}{2n^4} \pm\cdots\right] - \left[\frac{1}{n} + \frac{1}{2n^2}
- \frac{1}{3n^3} + \frac{1}{4n^4} \mp \cdots\right]
\\
&=\left[\frac{1}{2}-\frac{1}{3}\right]\frac{1}{n^3}-\left[\frac{1}{2}-\frac{1}{4}\right]\frac{1}{n^4}+\left[\frac{1}{2}-\frac{1}{5}\right]\frac{1}{n^5}-\left[\frac{1}{2}-\frac{1}{6}\right]\frac{1}{n^6}\mp\cdots\\
&= \frac{1}{6n^3} - \frac{1}{4n^4}+\frac{3}{10n^5}-\frac{1}{3n^6} \pm\cdots
\end{align*}
which is an alternating series whose terms decrease monotonically to
zero. A well-known theorem due to Leibniz states that if 
$$
S:=a_1-a_2+a_3-a_4\pm\cdots
$$
is an alternating series such that $a_n\geq 0$ and $a_n$ decreases monotonically to zero, then the series converges to a sum $S$ and if 
$$
S_n:=a_1-a_2+a_3-a_4\pm\cdots+(-1)^{n-1}a_n
$$
is the $n$-th partial sum, then the absolute value of the remainder $R_n$ satisfies:
$$
|R_n|:=|S-S_n|\leq a_{n+1}
$$ 
and the sign of $R_n$ is $(-1)^n$.
Therefore, by the Leibniz error estimate,
\begin{equation}
\label{eq:sigma-error} 
\frac{1}{6n^3} - \frac{1}{4n^4} < \sigma_n < \frac{1}{6n^3}\,.
\end{equation}
The standard estimate for the remainder from the integral
test is:  
$$
\int_{n+1}^\infty f(x) \,dx < R_n < \int_n^\infty f(x) \,dx
$$
where $R_n$ is the remainder
$$
R_n:=f(n+1)+f(n+2)+\cdots
$$
in the series $\sum_{n=1}^{\infty}f(n)$.
If we apply it to the series $\dsp \sum_{n=1}^\infty \frac{1}{6n^3}$ and 
$\dsp \sum_{n=1}^\infty \frac{1}{4n^4}$ we obtain
\begin{equation}
\label{eq:sigma-error-1} 
\frac{1}{12(n+1)^2} - \frac{1}{12n^3} < \sum_{n=1}^\infty \sigma_n
< \frac{1}{12n^2}\,.
\end{equation}
But,
\begin{align*}
\frac{1}{12(n+1)^2} - \frac{1}{12n^3}
&= \frac{1}{12n^2} - 2 \frac{1}{12n^3} + 3 \frac{1}{12n^4}-4 \frac{1}{12n^5}+5 \frac{1}{12n^6} \mp\cdots
- \frac{1}{12n^3}
\\
&= \frac{1}{12n^2} - 3 \frac{1}{12n^3} + 3 \frac{1}{12n^4}-4 \frac{1}{12n^5}+5 \frac{1}{12n^6} \mp\cdots\cdots
\\
&=\frac{1}{12n^2} -  \frac{1}{4n^3} + 3 \frac{1}{12n^4}-4 \frac{1}{12n^5}+5 \frac{1}{12n^6} \mp\cdots\cdots
\\
&> \frac{1}{12n^2} - \frac{1}{4n^3}
\end{align*}
since the series is alternating and the terms converge monotonically
to zero. Therefore, if we define
$$
\eps_n := \frac{1}{12n^2} -\sum_{n=1}^\infty\sigma_n
$$
we conclude that
\begin{equation}
\label{eq:epsilon} 
0 < \eps_n < \frac{1}{4n^3}
\end{equation}
as stated in the theorem. This completes the proof.


\section{Concluding remark}

Our method does not lead to an error term $O(1/n^4)$ since the terms
of order $1/n^3$ for $\sigma_n$ do not cancel. It would be desirable
to modify this geometric reasoning to achieve such a cancellation
(perhaps using telescopic cancellation if necessary).

\noindent Mark B. Villarino
\\
Escuela de Matemática, Universidad de Costa Rica,
\\
10101 San José, Costa Rica
\\
mark.villarino@ucr.ac.cr

\end{document}